\documentclass[12pt]{article}

\pdfoutput=1

\usepackage{amsmath,amsthm}
\usepackage{amssymb}
\usepackage{epsfig}

\usepackage{easybmat}
\usepackage{graphicx,hyperref}

\newtheorem{theorem}{Theorem}[section]
\newtheorem{lemma}[theorem]{Lemma}

\newtheorem{cor}[theorem]{Corollary}

\theoremstyle{definition}
\newtheorem{define}[theorem]{Definition}
\newtheorem*{remark}{Remark}

\newcommand{\pf}{\noindent {\bf Proof: }}
\newcommand{\enpf}{\hfill $\Box$ \vspace{.2in} }
\newcommand{\vs}{\vspace{0.15in}}

\begin{document}

\title{Stability of Relative Equilibria in the \\ Planar $n$-Vortex Problem}

\author{Gareth E. Roberts\thanks{
Dept. of Mathematics and Computer Science,
College of the Holy Cross, groberts@holycross.edu}}

\maketitle

\begin{abstract}
We study the linear and nonlinear stability of relative equilibria in the planar $n$-vortex problem, adapting the approach
of Moeckel from the corresponding problem in celestial mechanics.  After establishing some general theory,
a topological approach is taken to show that for the case of positive circulations, 
a relative equilibrium is linearly stable if and only if it is a nondegenerate minimum of the 
Hamiltonian restricted to a level surface of the angular impulse (moment of inertia).  
Using a criterion of Dirichlet's, this implies that any linearly stable relative equilibrium 
with positive vorticities is also nonlinearly stable.  Two symmetric families, the rhombus
and the isosceles trapezoid, are analyzed in detail, with stable solutions found in each case.
\end{abstract}

%{\bf MSC Classifications:} 70F10, 70F15, 37N05

\vspace{.1in}

{\bf Key Words:}  Relative equilibria, $n$-vortex problem, linear stability, nonlinear stability

%%%%%%%%%%%%%%%%%%%%%%%%%%%%%%%%%%%%%%%%%%%%%%%%%%%%%%%%%%%%%%%%%%%%%%%%%%%%%%%%%%%%
\section{Introduction}
%%%%%%%%%%%%%%%%%%%%%%%%%%%%%%%%%%%%%%%%%%%%%%%%%%%%%%%%%%%%%%%%%%%%%%%%%%%%%%%%%%%

In 2001, Kossin and Schubart~\cite{KossSchub} conducted numerical experiments describing
the evolution of thin annular rings with large vorticity as a model for the behavior seen in
the eyewall of intensifying hurricanes.   In a conservative, idealized setting, they find examples of 
``vortex crystals,'' formations of mesovortices that rigidly rotate as a solid body.
One particular formation of four vortices, situated very close to a rhombus configuration,
is observed to last for the final 18 hours of a 24-hour simulation (see Figure~4(a) of~\cite{KossSchub}).
Rigidly rotating polygonal configurations have also been found 
in the eyewalls of hurricanes in weather research and forecasting models
from the Hurricane Group at the National Center for Atmospheric Research
(see~\cite{davis} or the website~\cite{corb-web} for some revealing simulations).

It is natural to explore these rigidly rotating configurations in a dynamical systems setting by studying
relative equilibria of the planar $n$-vortex problem.  
Introduced by Helmholtz~\cite{helm} and later given a 
Hamiltonian formulation by Kirchhoff~\cite{kirchhoff},
the $n$-vortex problem is a widely used model for providing finite-dimensional approximations
to vorticity evolution in fluid dynamics~\cite{aref-int, newton}.   The general goal is to track the motion
of the vortices as points rather than focus on their internal structure and deformation, 
a concept analogous to the use of ``point masses'' in celestial mechanics.
A relative equilibrium is a special configuration of vortices that rigidly rotates about 
its center of vorticity.  In a rotating frame these solutions are fixed.
Given their persistence in many models of hurricanes, analyzing the stability of relative equilibria
in the planar $n$-vortex problem is of great interest and may have some practical significance.

Several researchers have studied the stability of relative equilibria, beginning
with the work of Lord Kelvin~\cite{kelvin1, kelvin2}, Gr\"{o}bli~\cite{grobli} and Synge~\cite{synge} 
on the well-known equilateral triangle solution (the analogue of Lagrange's solution in the three-body problem).
If $\Gamma_i$ represents the vorticity or circulation of the $i$-th vortex,  then placing the
vortices at the vertices of an equilateral triangle results in a relative equilibrium for any choice of
circulations (positive or negative).  This periodic solution is linearly stable if and only if the 
{\em total vortex angular momentum} $L$, defined as
$$
L \; = \;  \sum_{i < j}  \Gamma_i \Gamma_j ,
$$
is positive~\cite{synge}.  As we will show, the quantity~$L$ turns out to be very important in 
the analysis of stability.   We note that both Synge~\cite{synge} and Aref~\cite{aref-3vorts} show, using a special
coordinate system, that the equilateral triangle solution is actually nonlinearly stable whenever it is linearly
stable.  

Another interesting stability result concerns the regular $n$-gon ($n \geq 4$), a relative equilibrium only if all vortices have
the same strength.  It is linearly stable for $n \leq 6$, degenerate at $n=7$, and unstable for 
$n \geq 8$~\cite{thomson, have, aref-equil}.  By calculating the higher order terms in the normal form of the Hamiltonian,  
Cabral and Schmidt~\cite{cs, schm} show that in the special case $n=7$, 
the ``Thomson heptagon'' is locally Lyapunov stable.   It is possible to extend the regular $n$-gon relative 
equilibrium by adding an additional vortex of arbitrary circulation~$\nu$ at the center of the
ring.  Cabral and Schmidt show that, for all $n \geq 3$, this $1+n$-gon relative equilibrium is locally Lyapunov
stable as long as $\nu$ belongs to a bounded interval of values that depends on~$n$~\cite{cs}.
As was the case with the equilateral triangle, for each $n$, the interval where the $1+n$-gon is linearly stable coincides
with the one where the configuration is nonlinearly stable.   This same phenomenon where linear stability and
nonlinear stability coincide also occurs in many cases for point vortices on the sphere~\cite{boatt, csm, polz}.
Finally, we mention the recent work of Barry, Hall and Wayne considering the stability of relative equilibria containing one 
``dominant'' vortex ($\nu_0 >> 0$) and $n$ equal strength ($\Gamma_i = \epsilon$) but small vortices~\cite{bhw}.
Following the analogous work of Moeckel in the $n$-body problem~\cite{rick-dom}, they show that a relative
equilibrium with $\epsilon > 0$ ($\epsilon < 0$, resp.) in this setup is linearly stable provided the $n$ small 
vortices are situated at a nondegenerate local minimum (maximum, resp.) of a special potential function.

In this work, we follow Moeckel's~\cite{rick-stab} approach for studying linear stability of some symmetric configurations of relative
equilibria in the $n$-body problem.  Working in a rotating frame, the basic idea is to find invariant subspaces of the Hessian 
of the Newtonian potential (multiplied on the left by the inverse of a diagonal mass matrix), and use them to determine 
the linear stability of the relative equilibrium.   While this is somewhat complicated in the $n$-body case, 
it is straight-forward for the $n$-vortex problem, in part because the number of degrees of freedom is just~$n$, rather than~$2n$.  
Much of our analysis builds from a simple algebraic property of the logarithmic Hamiltonian of the planar $n$-vortex problem:
the Hessian anti-commutes with a key rotation matrix.

One of the interesting conjectures in the $n$-body problem, due to Moeckel and listed as Problem 16 in~\cite{albouy-pblms},
states that a linearly stable relative equilibrium must be a nondegenerate minimum of the Newtonian potential 
restricted to a level surface of the moment of inertia.  For the case of $n$ vortices with the same sign, we prove that
linear stability is in fact {\em equivalent} to being a nondegenerate minimum of the Hamiltonian $H$ restricted to a level surface
of the angular impulse (the analog of the moment of inertia).  Since the angular impulse is also a conserved quantity for
the $n$-vortex problem, it follows that, in the case of same signed circulations, linear stability implies nonlinear stability.
In other words, solutions that start near a linearly stable relative equilibrium stay close to it for all forward and backward time.
This helps explain many of the results described above where linear stability and nonlinear stability coincide.  
We also show that the collinear relative equilibria with same signed circulations are all unstable, and that for generic
choices of positive circulations, a linearly stable relative equilibrium always exists.
While previous works have focused on specific, often symmetric relative equilibria, our results are 
more general.

The paper is organized as follows.  In Section~2 we describe the $n$-vortex problem in rotating coordinates, highlight the key properties
of the Hamiltonian~$H$, and establish some general theory for determining the linear stability of a relative equilibrium.
The key quantity~$L$ is discussed and it is shown that relative equilibria with $L=0$ are always degenerate.  
Section~3 focuses on a topological approach to the problem and proves our main result connecting linear stability
and being a minimum of $H$ restricted to a level surface of the angular impulse.  The relation between linear stability
and nonlinear stability is carefully explained and proven.  Finally, in Section~4 we apply some of our theory to
analyze the linear stability of two families of symmetric relative equilibria in the four-vortex problem: a rhombus and an isosceles trapezoid.
In both cases, stable examples exists for positive circulations.  
The rhombus case is particularly interesting as it contains two one-parameter families of solutions, one of which goes through
a pitchfork bifurcation at a parameter value expressed as the root of a certain cubic polynomial.
Throughout this work, calculations and computations were confirmed symbolically using Maple~\cite{maple} and
numerically using Matlab~\cite{matlab}.

%%%%%%%%%%%%%%%%%%%%%%%%%%%%%%%%%%%%%%%%%
\section{Relative Equilibria}
\label{sec:rel-eq}
%%%%%%%%%%%%%%%%%%%%%%%%%%%%%%%%%%%%%%%%%

A system of $n$ planar point vortices with vortex strength $\Gamma_i \neq 0$ 
and positions $z_i \in \mathbb{R}^2$ evolves according to 
$$
\Gamma_i  \dot z_i \;  = \;  J \, \nabla_i H  \; = \;   J  \sum_{j \neq i}^n \frac{\Gamma_i \Gamma_j}{r_{ij}^2}(z_j - z_i), \quad 1 \leq  i \leq n
$$
where  
$$
H = -\sum_{i<j} \Gamma_i \Gamma_j \ln (r_{ij}), \quad r_{ij} = \|z_i - z_j\|
$$
is the Hamiltonian, $J$ is the standard $2 \times 2$ symplectic matrix
$
J=
\begin{bmatrix}
0 & 1\\
-1 & 0
\end{bmatrix}
$
and
$\nabla_i$ denotes the two-dimensional partial gradient with respect to $z_i$.
We will assume throughout that the  {\em total circulation}  $\Gamma = \sum_i  \Gamma_i $ is nonzero.  
The {\em center of vorticity} is then well-defined as $c =  \frac{1} {\Gamma} \sum_i \Gamma_i z_i$.
Without loss of generality, we take $c = 0$.  Another important quantity is 
the {\em total vortex angular momentum} $L$, defined earlier as
$
L  =   \sum_{i < j}  \Gamma_i \Gamma_j.
$
If all the vortex strengths are positive, then $L > 0$.  However, when vorticities of different
signs are chosen, $L \leq 0$ is possible.  The case $L = 0$ turns out to be of considerable
importance.

A configuration of vortices which rotates rigidly about its center of
vorticity is known as a relative equilibrium.
\begin{define}
A {\em relative equilibrium} is a solution of the form 
$$
z_i(t) = e^{- \omega J t} z_i(0) , \quad 1 \leq i \leq n,
$$ 
that is, a uniform rotation with angular velocity $\omega \neq 0$ 
around the origin.
\end{define}

The initial positions $z_i(0)$ of a relative equilibrium must satisfy 
\begin{equation}
- \omega \, \Gamma_i z_i(0) \; = \;  \nabla_i H \; = \;  \sum_{j \neq i}^n \frac{\Gamma_i \Gamma_j}{r_{ij}^2(0)} (z_j(0) - z_i(0) ), \quad \mbox{for each } i \in \{1, \ldots, n\}.
\label{eq:rel-equ}
\end{equation}
Let $z = (z_1, z_2, \ldots, z_n) \in \mathbb{R}^{2n}$ describe the vector of positions.
We will typically refer to the initial positions $z_0 = z(0) \in \mathbb{R}^{2n}$ as a relative equilibrium, 
rather than the full periodic orbit.

The important quantity $I = \frac{1}{2} \sum_{i=1}^n \, \Gamma_i \|z_i \|^2$ is the {\em angular impulse}, analogous to
the {\em moment of inertia} in the $n$-body problem.  
It is well known that $I$ is an integral of motion for the $n$-vortex problem~\cite{newton}.
The system of equations defined by~(\ref{eq:rel-equ}) can be written more
compactly as
\begin{equation}
\nabla H(z) + \omega \nabla I(z) \; = \; 0 .
\label{eq:cc}
\end{equation}
Thus, relative equilibria can be viewed as critical points of the Hamiltonian
restricted to a level surface of $I$, with $\omega$ acting as a Lagrange multiplier.  This provides a very useful topological 
approach to studying the problem.  

Note that if $z_0 \in \mathbb{R}^{2n}$ is a relative equilibrium, so
is $r z_0$ for any scalar $r$.  In this case, the angular velocity $\omega$ must be scaled by a 
factor of~$1/r^2$.  Moreover, if $e^{J \theta} z_0$ denotes $(e^{J \theta} z_1, \ldots , e^{J \theta} z_n)$, then
$e^{J \theta} z_0$ is also a relative equilibrium for any~$\theta$, with the same angular velocity.  Thus,
relative equilibria are not isolated.  When counting the number of solutions, it is customary to fix a scaling
(e.g., $I = 1$) and identify relative equilibria that are identical under rotation.  In this fashion,
one counts the number of equivalence classes of relative equilibria.  We also note that the stability
type of a particular relative equilibrium is unchanged by scaling or rotation.  In particular,
if $z_0$ is linearly stable, then so is any other relative equilibrium equivalent to $z_0$ under
scaling or rotation.

To analyze the stability of relative equilibria, we change to rotating coordinates.
Denote $M = {\rm diag}\{ \Gamma_1, \Gamma_1, \ldots, \Gamma_n, \Gamma_n \} $ as
the $2n \times 2n$ matrix of circulations.  Also, let $K$ be the $2n \times 2n$ block diagonal matrix
containing $J$ on the diagonal.  The equations of motion are thus written simply as
$
M \dot{z}  =  K \, \nabla H(z).
$
Under the transformation $z_i \mapsto e^{-\omega J t} z_i \; \forall i$,
the system is transformed into 
\begin{equation}
M \dot{z} \; = \;  K ( \nabla H(z) + \omega M z) ,
\label{eq:RotCoords}
\end{equation}
a uniformly rotating frame with period $2 \pi/\omega$.
As expected, a rest point of~(\ref{eq:RotCoords}) is a solution to
equation~(\ref{eq:cc}), a relative equilibrium.

\begin{remark}
Note that the system of differential equations in (\ref{eq:RotCoords}) is not written in the 
usual canonical form since the conjugate variables $x_i, y_i$ appear consecutively.
One could permute the variables and use the canonical $2n \times 2n$ matrix $J$ as opposed to $K$.
However, we choose this set-up since it is similar to the Newtonian $n$-body problem
and because a key property of the Hamiltonian is more readily apparent in these coordinates.
\end{remark}

The $\cdot$ in the following lemma represents the standard Euclidean dot product and
 $D^2 H(z)$ is the Hessian matrix of~$H$.
\begin{lemma}
The Hamiltonian $H$ has the following three properties:
\begin{enumerate}
\item[{\bf (i)}]
$ \nabla H(z) \cdot z  \; = \;  - L, $

\item[{\bf (ii)}]
$\nabla H(z) \cdot (K z) \; = \; 0$,

\item[{\bf (iii)}] 
$D^2H(z) \, K \; = \;  - K \, D^2H(z)$.

\end{enumerate}
\label{lemma:props}
\end{lemma}

\pf
Differentiating the identity $H(r z) = H(z) - L \ln|r|$ with respect to the scalar $r$ and evaluating at $r=1$ verifies property {\bf (i)}.
If $e^{Jt} z$ denotes $(e^{Jt} z_1, \ldots , e^{Jt} z_n)$, then the $SO(2)$ symmetry of $H$ gives
$H( e^{Jt} z ) =   H(z) \; \forall t$.  Differentiating this identity with respect to $t$ and evaluating at $t=0$ shows property {\bf (ii)}.
To prove the third property, notice that differentiation of
{\bf (i)} implies $-K D^2H(z) \, z =   K \, \nabla H(z)$, while
differentiation of {\bf (ii)} yields
$D^2H(z) K \, z =   K \, \nabla H(z)$.  Since $z$ is arbitrary, these two equations prove {\bf (iii)}.
\enpf

\vs

Using property {\bf (i)} and the fact that $I$ is homogeneous of degree~$2$, equation~(\ref{eq:cc}) implies that the angular
velocity of a relative equilibrium is $\omega = L/(2I)$.  If all the vortex strengths are positive,
then $\omega > 0$ is assured.  If $L = 0$, then we must also have $I = 0$ in order
to have a relative equilibrium.
Since $dI/dt  = -2 \nabla H(z) \cdot (K z)$, property {\bf (ii)} verifies that
the moment of inertia is a conserved quantity in the planar $n$-vortex problem.

The fact that the Hessian of the Hamiltonian and the rotation matrix $K$ anti-commute
will be profoundly important in the stability analysis of relative equilibria.  For positive vorticities,
it implies a complete factorization of the characteristic polynomial into even quadratic factors.
We note that property {\bf (iii)} holds true for any degree zero homogeneous potential with $SO(2)$ symmetry.

%%%%%%%%%%%%%%%%%%%%%%%%%%%%%%%%%%%%%%%%%%
\subsection{Linear Stability}
%%%%%%%%%%%%%%%%%%%%%%%%%%%%%%%%%%%%%%%%%%

We will assume throughout that $z_0 \in \mathbb{R}^{2n}$ represents a particular
relative equilibrium.  Linearizing system~(\ref{eq:RotCoords}) about $z_0$ gives the stability matrix
$$
B \; = \;  K( M^{-1} D^2 H(z_0) + \omega I ) ,
$$
where $I$ is the $2n \times 2n$ identity matrix.
Here we have used the fact that $K$ and $M^{-1}$ commute.
Through a permutation of the variables, our stability matrix $B$ is similar to the matrix used in~\cite{cs}.
The eigenvalues of $B$ determine the linear stability of the corresponding periodic solution.
Since the system is Hamiltonian, they come in pairs $\pm \lambda$.  To have 
stability, the eigenvalues must lie on the imaginary axis.  
The eigenvalues of $B$ can be described
simply in terms of the eigenvalues of $M^{-1} D^2 H(z_0)$.  The next two lemmas make this statement precise.

\begin{lemma}
The characteristic polynomials of $D^2H(z_0)$ and $M^{-1} D^2H(z_0)$ are even.  Moreover,
for each matrix, $v$ is an eigenvector with eigenvalue $\mu$ if and only if $Kv$ is an eigenvector
with eigenvalue $- \mu$.
\label{lemma:charpoly1}
\end{lemma}

\pf
We give the proof for $M^{-1} D^2H(z_0)$; the proof for $D^2H(z_0)$ is basically identical.
Note that ${\rm det}(K) = 1$ and $K^2 = -I$.  Let $q(\lambda) = {\rm det}(M^{-1} D^2H(z_0) - \lambda I)$ be the
characteristic polynomial for $M^{-1} D^2H(z_0)$.  Then, using property~{\bf (iii)} of Lemma~\ref{lemma:props}, we have 
\begin{eqnarray*}
q(\lambda) & = &  {\rm det}(K) \, {\rm det} (M^{-1} D^2H(z_0) - \lambda I) \, {\rm det}(K) \\
& = &  {\rm det}( \, K (M^{-1} D^2H(z_0) K + \lambda I )\\
& = &  {\rm det}(\, M^{-1} D^2H(z_0)  + \lambda I )\\
& = & q(-\lambda).
\end{eqnarray*}
We also have the following sequence of implications:
\begin{eqnarray*}
M^{-1} D^2H(z_0) \, v \; = \; \mu  v  &\Longleftrightarrow&  K M^{-1} D^2H(z_0) v \; = \; \mu K v \\
&\Longleftrightarrow&  -M^{-1} D^2H(z_0) K \, v \; = \; \mu K v \\
&\Longleftrightarrow&  M^{-1} D^2H(z_0) (K v) \; = \; - \mu (K v) , 
\end{eqnarray*}
which verifies the second part of the lemma.
\enpf

\begin{lemma}
Let $p(\lambda)$ denote the characteristic polynomial of the stability matrix $B$.
\begin{itemize}
\item[{\bf (a)}] Suppose that $v$ is a real eigenvector of  $M^{-1} D^2H(z_0)$ with eigenvalue $\mu$.  Then $\{ v, Kv \}$ is a real invariant subspace of $B$
and the restriction of $B$ to $\{ v, Kv \}$ is
\begin{equation}
\begin{bmatrix}
0 &\mu - \omega \\
\mu + \omega & 0 \\
\end{bmatrix}.
\label{eq:mat1}
\end{equation}
Consequently, $p(\lambda)$ has a quadratic factor of the form $\lambda^2 + \omega^2 - \mu^2$.

\item[{\bf (b)}]  Suppose that $v = v_1 + i \, v_2$ is a complex eigenvector of  $M^{-1} D^2H(z_0)$ 
with complex eigenvalue $\mu = \alpha + i \, \beta$.  Then $\{ v_1,v_2, Kv_1, Kv_2 \}$ is a real invariant subspace of $B$ 
and the restriction of $B$ to this space is
\begin{equation}
\begin{bmatrix}
0 & 0 & \alpha - \omega & \beta \\
0 & 0 & -\beta &  \alpha - \omega \\
\alpha + \omega & \beta & 0 & 0 \\
-\beta &\alpha + \omega & 0 & 0 
\end{bmatrix}.
\label{eq:mat2}
\end{equation}
Consequently, $p(\lambda)$ has a quartic factor of the form 
$(\lambda^2 + \omega^2 - \mu^2)(\lambda^2 + \omega^2 - \overline{\mu}^{\, 2})$.
\end{itemize}

\label{lemma:charpoly2}
\end{lemma}

\pf
First, suppose that $M^{-1} D^2H(z_0) v = \mu v$, with $v$ real.  Using Lemma~\ref{lemma:charpoly1}, this implies that
$B v = (\mu + \omega) K v$ and $B (Kv) =  (\mu - \omega) v$, verifying matrix~(\ref{eq:mat1}).  The characteristic polynomial
of matrix~(\ref{eq:mat1}) is $\lambda^2 + \omega^2 - \mu^2$ and consequently, this quadratic must be a factor of $p(\lambda)$.

Second, suppose that $M^{-1} D^2H(z_0) v = \mu v$, with $v = v_1 + i \, v_2$ and $\mu = \alpha + i \, \beta$.  
Using property~{\bf (iii)} of Lemma~\ref{lemma:props}, this implies that
$B v_1 = (\alpha + \omega) K v_1 - \beta K v_2$, $B v_2 = (\alpha + \omega) K v_2 + \beta K v_1$,
$B (Kv_1) =  (\alpha - \omega) v_1 - \beta v_2$ and $B (Kv_2) =  (\alpha - \omega) v_2 + \beta v_1$,
which verifies matrix~(\ref{eq:mat2}).  The characteristic polynomial
of matrix~(\ref{eq:mat2}) is $(\lambda^2 + \omega^2 - \mu^2)(\lambda^2 + \omega^2 - \overline{\mu}^{\, 2})$
and consequently, this quartic, which has real coefficients after expanding, must be a factor of $p(\lambda)$.
\enpf

Lemma~\ref{lemma:charpoly2} suggests a simple strategy for determining linear stability of a relative equilibrium~$z_0$.  
If enough eigenvectors of $M^{-1} D^2H(z_0)$ can be found, then the corresponding eigenvalues $\mu_i$
determine stability.  If $\mu_i \in \mathbb{R}$, then $|\mu_i| < |\omega|$ is required for linear stability.
If $\mu_i \in \mathbb{C}$, then Re$(\mu_i) = 0$ is required for linear stability.

Write the relative equilibrium as $z_0 = (z_1, z_2, \ldots, z_n) \in \mathbb{R}^{2n}$ and let $z_{ij} = (z_i - z_j)/r_{ij}$.
Direct computation shows that 
\begin{equation}
D^2H(z_0) \; = \;  
\begin{bmatrix}
A_{11} & A_{12} & \cdots & A_{1n} \\
\vdots &    &    &   \vdots \\
A_{n1} &  A_{n2}  & \cdots &  A_{nn}
\end{bmatrix}
\label{eq:D2H}
\end{equation}
where $A_{ij}$ is the $2 \times 2$ symmetric matrix 
\begin{equation}
A_{ij} \; = \;  \frac{ \Gamma_i \Gamma_j}{r_{ij}^2} [I - 2 z_{ij} z_{ij}^T] \quad \mbox{if } i \neq j,  \qquad
A_{ii} \; = \;  - \sum_{j \neq i} A_{ij}.
\label{eq:Aij}
\end{equation}
Note that $A_{ij} = A_{ji}$ and that for $i \neq j$,
$$
A_{ij} = 
\frac{\Gamma_i \Gamma_j}{r_{ij}^4} 
\begin{bmatrix}
(y_i - y_j)^2 - (x_i - x_j)^2 & -2(x_i - x_j)(y_i - y_j) \\[0.1in]
-2(x_i - x_j)(y_i - y_j) &  (x_i - x_j)^2 - (y_i - y_j)^2
\end{bmatrix} .
$$
The fact that $J$ commutes with each $A_{ij}$ gives another proof that $D^2H(z_0)$ and $K$ anti-commute.

Any relative equilibrium in the planar $n$-vortex problem will always have the four eigenvalues
$0, 0, \pm \omega i$.  We will refer to these eigenvalues as {\em trivial}.
They arise from the symmetry and integrals of the problem.
Due to conservation of the center of vorticity, the vectors $s = [1, 0, 1, 0, \ldots, 1, 0]$ and $Ks$ are in the kernel of
$D^2H(z_0)$.  This is easily verifiable from the structure of $D^2H(z_0)$.  
Consequently, $\mu = 0$ in Lemma~\ref{lemma:charpoly2} and $B$ has the eigenvalues
$\pm \omega i$.  In addition, due to the fact that relative equilibria are not isolated under rotation,
the vector $K z_0$ is in the kernel of $B$.  This can be verified analytically by differentiating
property {\bf (i)} of Lemma~\ref{lemma:props} and substituting in equation~(\ref{eq:cc}).  This gives the important identity
\begin{equation}
M^{-1} D^2H(z_0) \, z_0 \; = \;  \omega z_0,
\label{eq:evec-RE}
\end{equation}
which holds for any relative equilibrium~$z_0$.  Equation~(\ref{eq:evec-RE}) together with 
property {\bf (iii)} of Lemma~\ref{lemma:props} then shows that $BK z_0 = 0$, as expected.

For a given relative equilibrium $z_0$, let $V = \mbox{span} \{z_0, Kz_0 \}$.  This is an invariant subspace
for $B$ and the restriction of $B$ to $V$ is given by
$$
\begin{bmatrix}
0 & 0 \\
2 \omega & 0 \\
\end{bmatrix}.
$$
Thus, a relative equilibrium is always degenerate in our coordinates.
The off-diagonal term represents the fact that scaling $z_0$ gives another relative equilibrium but with a scaled angular velocity~$\omega$.  
We follow Moeckel's approach in~\cite{rick-stab} and define linear stability by restricting to a complementary subspace
of  $V$.

Recall that  $M = {\rm diag}\{ \Gamma_1, \Gamma_1, \ldots, \Gamma_n, \Gamma_n \} $.  It is important to note that, unlike
the $n$-body problem, $M$ will be indefinite if two vortices have circulations of opposite sign.
We say that the vectors $v$ and $w$ are {\em $M$-orthogonal} if $v^T M w = 0$.  
Next, let $V^\perp \subset \mathbb{R}^{2n}$ denote the $M$-orthogonal complement of $V$, that is,
$$
V^\perp \; = \;  \{ w \in  \mathbb{R}^{2n} : w^T M v = 0 \; \; \forall v \in V \} .
$$

\begin{lemma}
The vector space $V^\perp$ has dimension $2n - 2$ and is invariant under $B$.  
If $L \neq 0$, then $V \cap V^\perp = \{ 0 \}$.
\label{lemma:Vperp}
\end{lemma}

\pf
Note that $w \in V^\perp$ if and only if $w$ is in the null space of the $2 \times 2n$ matrix
\begin{equation}
\begin{bmatrix}
z_0^T M \\[0.07in]
z_0^T K M  \\
\end{bmatrix}.
\label{matrixVperp}
\end{equation}
Since the rank of this matrix is always two, the dimension of the kernel, and hence the dimension of $V^\perp$, is $2n - 2$.

Next, suppose that $w \in V^\perp$ and consider the vector $Bw$.  Using equation~(\ref{eq:evec-RE}) and property {\bf (iii)} of Lemma~\ref{lemma:props}, 
we have that  $B^T M z_0 = 0$ and $B^T M K z_0 = 2 \omega M z_0$.  This in turn implies that $(Bw)^T M z_0 = 0$
and $(Bw)^T M K z_0 = 2 \omega (w^T M z_0) = 0$.  Thus, $Bw \in V^\perp$ and $V^\perp$ is invariant under~$B$.

Finally, suppose that $v \in V \cap V^\perp$.  Write $v = a z_0 + b K z_0$ for some constants $a$ and $b$.  Since $v$ is also in $V^\perp$, we have
\begin{eqnarray}
(a z_0 + b Kz_0)^T M z_0   & = &  0, \label{eq:perp1} \\
(a z_0 + b Kz_0)^T M K z_0   & = &  0 . \label{eq:perp2} 
\end{eqnarray}
But $z^T M K z = 0$ for any $z \in \mathbb{R}^{2n}$ 
and since $z_0$ is a relative equilibrium, $z_0^T M z_0 = 2I = L/\omega$.  Thus, equation~(\ref{eq:perp1}) reduces
to $aL = 0$ and equation~(\ref{eq:perp2}) becomes $bL = 0$.  Since $L \neq 0$ was assumed, we have $a = b = 0$ and $v = 0$, as desired.
\enpf

As long as $L \neq 0$,
Lemma~\ref{lemma:Vperp} allows us to define linear stability with respect to the $M$-orthogonal complement of the  subspace $V$.
Instead of working on a reduced phase space (eliminating the rotational symmetry), we will stay in the
full space and define linear stability by restricting~$B$ to~$V^\perp$.

\begin{define}
A relative equilibrium $z_0$ always has the four trivial eigenvalues $0, 0, \pm \omega i$.  We call $z_0$ {\em nondegenerate}
if the remaining $2n-4$ eigenvalues are nonzero.  A nondegenerate relative equilibrium  is {\em spectrally stable} if the nontrivial 
eigenvalues are pure imaginary, and {\em linearly stable} if, in addition, the restriction of the stability matrix $B$ to $V^\perp$ 
has a block-diagonal Jordan form with blocks
$
\begin{bmatrix}
0 & \beta_i \\
-\beta_i & 0 \\
\end{bmatrix}.
$
\end{define}

%%%%%%%%%%%%%%%%%%%%%%%%
\subsection{The Special Case $L = 0$}
%%%%%%%%%%%%%%%%%%%%%%%%

Under this definition, a relative equilibrium ($\omega \neq 0$) with $L = 0$ will {\em always} be degenerate.  
This fact is interesting in its own right; changes in stability and bifurcations are expected to
occur when $L = 0$.  This is precisely what transpires for the equilateral triangle solution
as well as for the rhombus B family (see Section~\ref{sec:rhombus}).

\begin{theorem}
Suppose that $z_0$ is a relative equilibrium for a choice of circulations with $L = 0$.  Then two of the nontrivial
eigenvalues for $z_0$ are zero and $z_0$ is degenerate.
\label{thm:degen}
\end{theorem}

\pf
Suppose that $L = 0$.   By Lemma~\ref{lemma:Vperp},
$V^\perp$ has dimension $2n - 2$ and is invariant under $B$.  
If $z_0$ satisfies equation~(\ref{eq:cc}) for some $\omega \neq 0$, then we must have
$ 2I = z_0^T M z_0 = 0$.  This leads to the rather perverse setting where $V = \mbox{span} \{z_0, Kz_0 \}$ is
contained in $V^\perp$, its $M$-orthogonal complement.  Using matrix~(\ref{matrixVperp}) as a guide, consider 
a basis $\xi$ of $\mathbb{R}^{2n}$ of the form 
$$
\xi \; = \;  \{ Mz_0, MKz_0, \, z_0, Kz_0, s, Ks, w_1, \ldots , w_{2n-6} \}
$$
where the $w_i$ represent a set of $2n - 6$ linearly independent vectors that complete
a basis for $V^\perp$.  Note that the first four vectors in $\xi$ are mutually orthogonal.

To find the stability matrix $B$ written with respect to $\xi$, write
\begin{eqnarray}
B M z_0 & = &  a_1 Mz_0 + a_2 MKz_0 + a_3 z_0 + a_4 Kz_0 + \cdots  , \label{eq:BMz} \\
B MKz_0 & = & b_1 Mz_0 + b_2 MKz_0 + b_3 z_0 + b_4 Kz_0 + \cdots . \label{eq:BMKz} 
\end{eqnarray}
Taking the standard dot product of both sides of equation~(\ref{eq:BMz}) with the vectors $Mz_0$ and $MKz_0$,
respectively, gives $0 = a_1 z_0^T M^2 z_0$ and $2 \omega  z_0^T M^2 z_0 = a_2 z_0^T M^2 z_0$, respectively.  
Since  $z_0^T M^2 z_0 > 0$, we have $a_1 = 0$ and $a_2 = 2 \omega$.   Similarly,
taking the dot product of both sides of equation~(\ref{eq:BMKz}) with the vectors $Mz_0$ and $MKz_0$,
respectively, gives $0 = b_1 z_0^T M^2 z_0$ and $0 = b_2 z_0^T M^2 z_0$, respectively.  
It follows that $b_1 = b_2 =0$.  Thus, the stability matrix $B$ written with respect to the basis~$\xi$ has the form

$$
\left[
\begin{BMAT}(c){c:c:c}{c:c:c}
\begin{BMAT}(rc){cc}{cc}
0 & 0 \\
2 \omega & 0 
\end{BMAT}  & 0 & 0\\
\ast &
\begin{BMAT}(rc){cc}{cc}
\begin{BMAT}(rc){cc}{cc}
0 & 0 \\
2 \omega & 0 
\end{BMAT}  & 0 \\
0 &  \begin{BMAT}(rc){cc}{cc}
0 & - \omega \\
\omega & 0 
\end{BMAT} 
\end{BMAT} & \ast \\
\ast & 0 & \ast
\end{BMAT}
\right],
$$
and $z_0$ has two additional zero eigenvalues in addition to the expected eigenvalues $0, 0, \pm \omega i$.
\enpf

%%%%%%%%%%%%%%%%%%%%%%%%%%%%%%%%%%%%
\section{A Topological Approach}
%%%%%%%%%%%%%%%%%%%%%%%%%%%%%%%%%%%%

In this section we present our main results, focusing specifically on the case $\Gamma_j > 0 \; \forall j$.  
Recall that the angular velocity $\omega = L/(2I)$, so $\omega > 0$ is assured for this case.  We also have that $M$ is positive definite.
Thus we can define an inner product by
$$
<v,w> \; = \;  v^T M w .
$$
Let $\delta_j = [[\frac{j+1}{2}]]$, where $[[ \cdot ]]$ is the greatest integer function.  Denote $\{e_j\}$ as the standard basis vectors of $\mathbb{R}^{2n}$ 
and let  $p_j = (1/\sqrt{ \Gamma_{\delta_j}}) e_j$.  The vectors $\{p_1, p_2, \ldots, p_{2n} \}$ form an $M$-orthonormal basis
of $\mathbb{R}^{2n}$, that is, $<p_i,p_j> = \delta_{ij}$, where $\delta_{ij}$ is the usual Kronecker delta function.
The key matrix $M^{-1} D^2H(z_0)$ is symmetric with respect to this basis.  Consequently, all its eigenvalues
are real and it has a full set of $M$-orthogonal eigenvectors.  By Lemma~\ref{lemma:charpoly1}, the eigenvalues
are of the form $\pm \mu_j$.  It follows that the characteristic polynomial of the stability matrix $B$
factors completely into quadratic factors of the form $\lambda^2 + \omega^2 - \mu_j^2$.

\begin{theorem}
If $\Gamma_j > 0 \; \forall j$, then a relative equilibrium $z_0$ is linearly stable if and only if 
$| \mu_j | < \omega$ for all eigenvalues $\mu_j$ of $M^{-1} D^2H(z_0)$.
\label{thm:linstab}
\end{theorem}

\pf
For the case $\Gamma_j > 0 \; \forall j$, all the eigenvalues of $M^{-1} D^2H(z_0)$ are real and come in pairs
$\pm \mu_j$.  Part {\bf (a)} of Lemma~\ref{lemma:charpoly2} now applies repeatedly and the characteristic polynomial of $B$ factors as
$$
p(\lambda) \; = \;  \lambda^2 (\lambda^2 + \omega^2) \prod_{j = 1}^{n-2} (\lambda^2 + \omega^2 - \mu_j^2).
$$
If $z_0$ is linearly stable, then $\omega^2 - \mu_j^2 > 0$ is required for all~$j$.  Conversely, if $| \mu_j | < \omega$
for all~$j$, we are assured that all the nontrivial eigenvalues are nonzero and lie on the imaginary
axis.  Thus, $z_0$ is nondegenerate and since $M^{-1} D^2H(z_0)$ has a full set of eigenvectors of the form $\{v_j, Kv_j \}$,
the restriction of the stability matrix $B$ to $V^\perp$ has a block-diagonal Jordan form with blocks
$$
\begin{bmatrix}
0 & \sqrt{ \omega^2 - \mu_j^2}  \\
-\sqrt{ \omega^2 - \mu_j^2}   & 0 \\
\end{bmatrix}.
$$
Consequently, $z_0$ is linearly stable.
\enpf

\begin{remark}
\begin{enumerate}
\item
For the case of positive circulations, if $z_0$ is linearly stable, then the angular velocity of each component of the linearized 
flow is less than or equal to the angular velocity of $z_0$.
If we assume further that $\mu_j \neq 0 \; \forall j$, then, working on a reduced space 
that eliminates the $SO(2)$ symmetry, the Lyapunov center theorem applies to produce a one-parameter 
family of periodic orbits emanating from~$z_0$~\cite{meyer}.

\item  Theorem~\ref{thm:linstab} also applies to the case $\Gamma_j < 0 \; \forall j$, although here we have $\omega < 0$ so the condition for
linear stability becomes $| \mu_j | < - \omega$.

\item  Due to the structure of the Jordan form, it is not possible to have Krein collisions for this case.  In other words, stability can only
be lost by having a pair of pure imaginary eigenvalues merge together at 0 and then move onto the real axis; eigenvalue pairs that collide
on the imaginary axis will pass through each other.   With positive
circulations, the linear stability and spectral stability of $z_0$ are actually equivalent concepts.

\end{enumerate}
\end{remark}

Next we follow a topological approach to prove our main result.
Recall that $I = \frac{1}{2} \sum_{i=1}^n \, \Gamma_i \|z_i \|^2$ and that a relative equilibrium 
is a critical point of $H$ restricted to the level surface $I = I_0$.  Assuming all the
vortex circulations are positive, $I = I_0$ is an ellipsoid or topological sphere in $\mathbb{R}^{2n}$.
If we restrict $H$ to this sphere, it turns out that the local (nondegenerate) minima of $H$ are exactly the
linearly stable relative equilibria.  

First we recall some important definitions.  Suppose that $z_0$ is a relative equilibrium 
with $I(z_0) = I_0$.  Let $S$ denote the surface $I = I_0$ and let $\gamma \equiv \gamma(\theta)$ be the closed curve $e^{J \theta} z_0$
of relative equilibria.  The curve~$\gamma$ (a circle) lies on $S$ and $H$ is constant on~$\gamma$.
The tangent space of $S$ at $z_0$, denoted $T_{z_0}S$, is simply the $M$-orthogonal complement
of $z_0$.  In other words, $T_{z_0}S = \{v \in \mathbb{R}^{2n}:  v^T M z_0 = 0 \}$.  The tangent vector
to $\gamma$ at $z_0$ is just $K z_0$ and since $<Kz_0, z_0> = 0$, it lies in $T_{z_0}S$, as expected.

Define the Lagrangian function $G(z) = H(z) + \omega I(z)$.  By equation~(\ref{eq:cc}), $z_0$ is
a critical point of~$G$.  The type of critical point is determined by examining the quadratic form
$q_{z_0}(v)$ associated to the Hessian matrix of~$G$~\cite{spring}.  To work with this quadratic form,
it is easier to remain in $\mathbb{R}^{2n}$ rather than use local coordinates.  We have
$$
q_{z_0}(v) \; = \;  v^T D^2G(z_0) v \; = \;  v^T(D^2H(z_0) + \omega M) v, \quad v \in T_{z_0}S.
$$
The {\em nullity} of $q_{z_0}$ is the dimension of the null space of $D^2G$ and
the {\em index} (or Morse index) of $q_{z_0}$ is the dimension of the maximal subspace
for which $q_{z_0}(v) < 0$.   The connection between the stability matrix $B$ and the quadratic form
$q_{z_0}$ is given by the equation $-MKB = D^2G$.  It follows that $z_0$ is a nondegenerate 
relative equilibrium if and only if the nullity of $q_{z_0}$ is one.   Therefore, treating $z_0$ as a critical
point of $H$ restricted to $S$, we will call $z_0$ a nondegenerate
critical point if the nullity of $q_{z_0}$ is one.  In this case, since $H$ is constant along $\gamma$, if the index of $q_{z_0}$ is zero, 
then $z_0$ is a nondegenerate minimum of $H$ restricted to $S$.

\begin{theorem}
If $\Gamma_j > 0 \; \forall j$, then a relative equilibrium $z_0$ is linearly stable if and only if 
it is a nondegenerate minimum of $H$ subject to the constraint $I = I_0$.
\label{thm:main}
\end{theorem}

\pf
Suppose that $z_0$ is linearly stable.  Since $\Gamma_j > 0 \; \forall j$, the matrix $M^{-1} D^2H(z_0)$ has an
$M$-orthonormal set of eigenvectors of the form $\{v_j, Kv_j \}$ with eigenvalues $\pm \mu_j$.
Let us choose $v_1 = (1/\sqrt{2 I_0} \, ) z_0$ and $v_2 = (1/\sqrt{\Gamma} \,) [1, 0, \ldots , 1, 0]^T$.  
The vectors $\{Kv_1, v_2, Kv_2, \ldots , v_n, Kv_n\}$ form an $M$-orthonormal basis for $T_{z_0}S$.  
Note that 
$$
D^2G(z_0) v_j \; = \;  (\mu_j + \omega) M v_j  \qquad  \mbox{and } \qquad  D^2G(z_0) K v_j \; = \; (-\mu_j + \omega) M K v_j.
$$

Given any vector $v \in T_{z_0}S$, write it as $v = \sum_{j=2}^n c_j v_j + \sum_{j=1}^n d_j (K v_j)$.  We then compute that
$$
q_{z_0}(v) \; = \;  \omega (c_2^2 + d_2^2) + \sum_{j=3}^n [ c_j^2 (\mu_j + \omega) + d_j^2 (-\mu_j + \omega)  ].
$$
By Theorem~\ref{thm:linstab}, we have $- \omega < \mu_j < \omega \; \forall j$.  It follows that
$q_{z_0}(v) \geq 0$ for all $v \in T_{z_0} S$.  Moreover, $q_{z_0}(v) = 0$ if and only if
$v \in \mbox{span} \{K z_0 \}$, a direction tangent to the curve of critical points $\gamma$.
Thus, the nullity of $q_{z_0}$ is one and the index is zero. Hence, $z_0$ 
is a nondegenerate minimum of $H$ subject to the constraint $I = I_0$.

Conversely, suppose that $z_0$ is a nondegenerate minimum of $H$ restricted to $S$ and by contradiction,
suppose that $z_0$ is not linearly stable.  From Theorem~\ref{thm:linstab}, there must exist an eigenvector
$v$ of $M^{-1} D^2H(z_0)$ with eigenvalue $\mu$ such that $v \not \in \mbox{span} \{ z_0, K z_0 \}$
and $\mu \geq \omega$.   (If $\mu < 0$, then replace $v$ with $Kv$.)  We can also scale $v$ so that
$v^T M v = 1$.  We then have $Kv \in T_{z_0}S$ and 
$q_{z_0}(Kv) = \omega - \mu \leq 0$.  But this contradicts the fact that 
$z_0$ is a nondegenerate minimum of $H$ subject to the constraint $I = I_0$.
\enpf

\begin{remark}
As mentioned in the introduction, the discovery of this theorem was motivated by the related conjecture of Moeckel's in the $n$-body problem
(Problem 16 in~\cite{albouy-pblms}).   Unlike the case of $n$~vortices, the converse direction is clearly false for the
$n$-body problem.  For example, Lagrange's equilateral triangle solution in the three-body problem is a minimum for any choice of masses 
but is linearly stable only for a small set of mass values.
\end{remark}

For positive circulations, the existence of a minimum of $H$ restricted to $I = I_0$ is clear.  As $r_{ij} \rightarrow 0$,
$H \rightarrow \infty$ and since $S$ is compact, $H$ restricted to $I = I_0$ is bounded below.  It follows that $H$ restricted to $S$
has a global minimum, and therefore a relative equilibrium $z_0$ exists.  However, for special choices
of the circulations, $z_0$ may be degenerate.  For example, the relative equilibrium consisting of a regular $n$-gon (equal strengths) 
and a central vortex of circulation $\Gamma_0$ is degenerate for certain values of $\Gamma_0$.  The number of degenerate
cases grows with~$n$ (see~\cite{cs}).  Thus, the existence of a linearly stable relative equilibrium is assured for 
most choices of circulations.  It turns out that this stable relative equilibrium is never a collinear
configuration.

Given $n$ positive circulations, there exist exactly $n!/2$ collinear relative equilibria, one for each ordering
of the vortices (the factor of $1/2$ occurs because configurations identical under a $180^\circ$ rotation are identified).
This follows by generalizing a well-known result of Moulton's for the $n$-body problem~\cite{moulton, oneil}.
It can be shown that each collinear relative equilibrium is a minimum of $H$ restricted to the 
sub-manifold of collinear configurations on $S$.  However, in directions normal to the collinear manifold,
the value of $H$ decreases.  In other words, collinear relative equilibria are minima restricted to the
collinear case but saddles in the full planar problem.  Intuitively, as the vortices move off the line,
the mutual distances increase and the value of $H$ decreases.

A rigorous proof of this fact is subtle.  In~\cite{palmore}, Palmore states that 
the collinear relative equilibria have index~$n-2$ in the planar $n$-vortex problem, 
although no proof is given.  The same result is true for the planar $n$-body problem and 
follows from a very clever argument due to Conley, articulated by Pacella in~\cite{pacella}.
It is straight-forward to check that the proof given in~\cite{pacella} generalizes to the 
planar $n$-vortex problem.  The basic idea is to view the linear flow corresponding to 
the normal directions as a flow on the space of lines through the origin.  One can find a
``cone'' about the trivial eigenvector $Kz_0$ that is mapped inside itself.  This shows that
the eigenvalues of $M^{-1} D^2H(z_0)$ corresponding to the normal directions are
all less than $- \omega$, which means the index of $q_{z_0}$ is $n-2$.
The reader is encouraged to see~\cite{pacella} or the lecture notes of Moeckel~\cite{rick:notes} 
for more details.  Since collinear relative equilibria with positive vorticities are never local minima, we have the following result.

\begin{cor}
Any collinear relative equilibrium with positive vorticities is unstable.  All the nontrivial eigenvalues come in
real pairs $\pm \alpha_j$.
\label{cor:collinear}
\end{cor}

The study of collinear relative equilibria with circulations of mixed sign is far more subtle, 
as discovered by Aref in his comprehensive study of three-vortex relative equilibria~\cite{aref-stab3}.  
In this case, it is possible to have linearly stable solutions.  We note that Aref's result for the case of
three collinear vortices with positive circulations is in agreement with Corollary~\ref{cor:collinear}.
Since collinear relative equilibria are saddles, there must exist a strictly planar global minimum
for~$H$ restricted to~$S$.  This yields another nice corollary of Theorem~\ref{thm:main}.

\begin{cor}
For a generic choice of positive circulations, there exists a non-collinear,
linearly stable relative equilibrium.
\label{cor:exist}
\end{cor}

\begin{remark}
Corollary~\ref{cor:exist} is known to be false for the $n$-body problem.  In fact, if $n \geq 24,306$, then {\em all} of the equal-mass relative equilibria
are unstable~\cite{g:equmass}.
\end{remark}

Since both $H$ and $I$ are conserved along solutions of the planar $n$-vortex problem, Theorem~\ref{thm:main}
actually has deeper implications involving nonlinear stability.  Suppose that
$z_0$ is a linearly stable relative equilibrium.  Let $\gamma = e^{- \omega J t} z_0$ be the periodic solution 
corresponding to $z_0$.  Each point on $\gamma$ satisfies equation~(\ref{eq:cc}) and is locally
a nondegenerate minimum of $H$ restricted to $I = I_0 = I(z_0)$.
Modifying a classical result of Dirichlet's~\cite{dirichlet}, we will show that solutions starting close to $z_0$
will stay close to $\gamma$ (as a set) for all time (forward and backward).

To be precise, let $\phi (t, \zeta_0)$ denote
the solution to the differential equation~(\ref{eq:RotCoords}) with initial condition~$\zeta_0$.
Define the distance between a point $\zeta_0$ and a compact set $Z$ to be
$$
d(\zeta_0, Z) \; = \; \mbox{min} \, \{ ||\zeta_0 - z|| : z \in Z \}.
$$
We will call the relative equilibrium $z_0$ {\em nonlinearly stable} or simply {\em stable} 
if for any $\epsilon > 0$, there exists a $\delta > 0$ such that
$||\zeta_0 - z_0|| < \delta$ implies $d(\phi(t, \zeta_0), \gamma) < \epsilon$ for all $t \in \mathbb{R}$.
Since equation~(\ref{eq:RotCoords}) is the planar $n$-vortex problem in a rotating frame,
solutions starting in a neighborhood of a stable relative equilibrium 
will stay close in the original coordinates as well.

\begin{theorem}
Suppose $\Gamma_i > 0 \; \forall i$.  Then any linearly stable relative equilibrium is also nonlinearly stable.
\label{thm:stable}
\end{theorem}

\pf
We follow the related proof given in~\cite{meyer}.  Suppose that $z_0$ is a linearly stable relative
equilibrium and define $\gamma$ as above.  
Fix $\epsilon > 0$.    Using property~{\bf (ii)} of Lemma~\ref{lemma:props}, it is straight-forward to check
that $G = H + \omega I$ is a conserved quantity for system~(\ref{eq:RotCoords}).  
By adding a constant to the Hamiltonian $H$, we can assume, without loss of generality,
that $G(z_0) = 0$.   Since $z_0$ is linearly stable, Theorem~\ref{thm:main} implies that $G$
is strictly positive in a neighborhood of $z_0$ except along $\gamma$ where it vanishes.   The same is true for any point on $\gamma$.
Let $Y_\eta = \{ z \in \mathbb{R}^{2n}: 0 < d(z, \gamma) < \eta \}$ be the open tube of radius $\eta$
surrounding~$\gamma$.  Since $\gamma$ is compact, there exists an $\eta > 0$ such that $G$ is positive on $Y_\eta$.

Let $\epsilon^\ast = \mbox{min}\{\epsilon, \eta \}$ and let $c^\ast = \mbox{min}\{ G(z): z \in \partial Y_{\epsilon^\ast}  \}$.   Thus, $c^\ast > 0$.
Since $G$ is continuous away from collisions, there exists a $\delta$-neighborhood of $z_0$ 
where $G <  c^\ast$.  For any $\zeta_0$ in this neighborhood, we have $G(\phi(t,\zeta_0)) <  c^\ast$ for all $t$.  
If there exists some time $t^\ast$ where $d(\phi(t^\ast,\zeta_0), \gamma) = \epsilon^\ast$, then $G(\phi(t^\ast,\zeta_0)) \geq c^\ast$, which is a contradiction.
Therefore, we have $d(\phi(t,\zeta_0), \gamma) < \epsilon^\ast \leq \epsilon$ for all $t$, as desired.
\enpf

%%%%%%%%%%%%%%%%%%%%%%%%%%%%%%%%%%%%%%%%%
\section{Some Symmetric Examples}
%%%%%%%%%%%%%%%%%%%%%%%%%%%%%%%%%%%%%%%%%

In this section we apply some of the general theory established above to some specific
examples of relative equilibria in the planar $n$-vortex problem.  The examples
considered are all symmetric, which allows for a simpler analysis of the stability.

%%%%%%%%%%%%%%%%%%%%%%%%%%%%%%%%%%%%%%%%%
\subsection{The Equilateral Triangle}
%%%%%%%%%%%%%%%%%%%%%%%%%%%%%%%%%%%%%%%%%

We begin with the well-known equilateral triangle solution in the three-vortex problem.
Placing three vortices of any strength at the vertices of an equilateral triangle yields
a relative equilibrium.   According to Newton~\cite{newton}, Lord Kelvin was the first to conclude that the equilateral
triangle solution was stable for the case of identical circulations~\cite{kelvin1, kelvin2}.
The general case was handled later by Synge~\cite{synge} who used trilinear coordinates
to show that it was linearly stable if and only if $L > 0$.  We reproduce this result here using our setup and
applying Lemma~\ref{lemma:charpoly2}.

Given any three circulations $\Gamma_1, \Gamma_2, \Gamma_3$,
let $\hat{z_1} = (1,0), \hat{z_2} = (-1/2,\sqrt{3}/2), \hat{z_3} = (-1/2, -\sqrt{3}/2)$
and $\hat{c} = \frac{1}{\Gamma} \sum_i  \Gamma_i \hat{z_i}$.
Then the vortex positions $z_i = \hat{z_i} - \hat{c}$ will have center of vorticity at the
origin and $z_0 = (z_1, z_2, z_3)$ is a relative equilibrium with angular velocity
$\omega = \Gamma/3$.  Using equations (\ref{eq:D2H}) and~(\ref{eq:Aij}),
we compute that $M^{-1} D^2H(z_0)$ is given by
$$
\frac{1}{6}
\begin{bmatrix}
\Gamma_2 + \Gamma_3  & \sqrt{3} \, (\Gamma_3 - \Gamma_2)  & - \Gamma_2 & \sqrt{3} \, \Gamma_2  &  - \Gamma_3  &  - \sqrt{3} \, \Gamma_3 \\[0.05in]
\sqrt{3} \, (\Gamma_3 - \Gamma_2)  & - (\Gamma_2 + \Gamma_3) & \sqrt{3} \, \Gamma_2  & \Gamma_2  & - \sqrt{3} \, \Gamma_3  &  \Gamma_3 \\[0.05in]
- \Gamma_1 &   \sqrt{3} \, \Gamma_1  &  \Gamma_1 - 2 \Gamma_3  &  - \sqrt{3} \, \Gamma_1  &  2 \Gamma_3  & 0 \\[0.05in]
\sqrt{3} \, \Gamma_1 & \Gamma_1  &   - \sqrt{3}\,  \Gamma_1  &  - ( \Gamma_1 - 2 \Gamma_3)  &  0  & -2 \Gamma_3 \\[0.05in]
- \Gamma_1 &   -\sqrt{3} \, \Gamma_1  &   2 \Gamma_2  & 0 &  \Gamma_1 - 2 \Gamma_2  &   \sqrt{3} \, \Gamma_1 \\[0.05in]
- \sqrt{3} \, \Gamma_1 & \Gamma_1    &  0  & -2 \Gamma_2 &  \sqrt{3} \,  \Gamma_1  &  - ( \Gamma_1 - 2 \Gamma_2)  
\end{bmatrix} .
$$
We know that $M^{-1} D^2H(z_0)$ has the four trivial eigenvalues $0, 0, \omega, -\omega$.  Let $\pm \mu$ be the remaining
two eigenvalues, where $\mu$ may be real or pure imaginary.  Instead of searching for an eigenvector, we follow Moeckel's approach in~\cite{rick-stab} and use
a formula involving the trace of the square of $M^{-1} D^2H(z_0)$ to obtain the remaining eigenvalues.
The characteristic polynomial of $M^{-1} D^2H(z_0)$ is 
\begin{equation}
\lambda^6 + a \lambda^4 + b \lambda^2 \; = \;  \lambda^6 - (\mu^2 + \omega^2) \lambda^4 + \mu^2 \omega^2 \lambda^2.
\label{eq:tri-poly}
\end{equation}
On the other hand, the coefficient $a$ in~(\ref{eq:tri-poly}) can be found using
$a = \frac{1}{2}[ (\mbox{Tr}(C) )^2 - \mbox{Tr}(C^2) ]$ where $C = M^{-1} D^2H(z_0)$ and
$\mbox{Tr}(\ast)$ is the trace of a matrix.  This formula applies nicely here because $\mbox{Tr}(C)=0$.
A straight-forward computation gives $\mbox{Tr}(C^2) = \frac{4}{9} \Gamma^2 - \frac{2}{3} L$.
Then we have
$$
 - \frac{2 \Gamma^2}{9} + \frac{L}{3}  \; = \; a \; = \;  - (\mu^2 + \omega^2) \; = \;  -\mu^2 - \frac{ \Gamma^2}{9}
\qquad  \Longrightarrow \qquad
\mu^2 \; = \;  \frac{\Gamma^2}{9} - \frac{L}{3} .
$$
Expanding $\mu^2$ in terms of the circulations $\Gamma_j$ yields
$$
\mu^2 \; = \;  \frac{1}{9} \left(  \Gamma_1^2 + \Gamma_2^2 + \Gamma_3^2 - \Gamma_1 \Gamma_2 - \Gamma_1 \Gamma_3 - \Gamma_2 \Gamma_3 \right)
\; = \;  \frac{1}{36} \left[  (\Gamma_1 + \Gamma_2 - 2 \Gamma_3)^2 + 3(\Gamma_1 - \Gamma_2)^2 \right],
$$
which shows that $\mu^2 \geq 0$ with equality if and only if the three circulations are equal.
Note that $\omega^2 - \mu^2 = L/3$.  Applying Lemma~\ref{lemma:charpoly2}, we see that the nontrivial eigenvalues of $z_0$
are $\pm \sqrt{ -L/3 }$.  Therefore, the equilateral triangle is linearly stable when the circulations satisfy
$L > 0$,  unstable when $L < 0$ and degenerate when $L = 0$.  This agrees with the result of Synge~\cite{synge}.
We also point out that for the case of equal circulations, $\mu = 0$ and the nontrivial eigenvalues are $\pm \omega i$.
In this case, the Floquet multipliers for the periodic solution (in the full phase space) are all equal to one.

%%%%%%%%%%%%%%%%%%%%%%%%%%%%%%%%%%%%%%%%%
\subsection{The Rhombus Families}
\label{sec:rhombus}
%%%%%%%%%%%%%%%%%%%%%%%%%%%%%%%%%%%%%%%%%

For the planar four-vortex problem, Hampton, Santoprete and the author showed that
there exists two families of relative equilibria where the configuration is a rhombus~\cite{HRS}.   
Set  $\Gamma_1 = \Gamma_2 = 1$ and $\Gamma_3 = \Gamma_4 = m$, where
$m \in (-1,1]$ is a parameter.    Note that $L = m^2 + 4m + 1$ vanishes when 
$m = -2 + \sqrt{3} \approx -0.268$, so this is an expected bifurcation value.
Place the vortices at $z_1 = (1,0), z_2 = (-1,0),
z_3 = (0,y)$ and $z_4 = (0,-y)$, forming a rhombus with diagonals lying on the
coordinate axis.  This configuration is a relative equilibrium provided that
\begin{equation}
y^2 \; = \;  \frac{1}{2} \left( \beta \pm \sqrt{\beta^2 + 4m} \, \right), \quad \beta \; = \;  3(1-m) ,
\label{eq:rhom-my}
\end{equation}
or equivalently, 
\begin{equation}
m \; = \;  \frac{3y^2 - y^4}{3y^2 - 1} .
\label{eq:rhom-myA}
\end{equation}
The angular velocity is given by
$$
\omega \; = \;  \frac{m^2 + 4m + 1}{2(1+my^2)} \; = \;  \frac{1}{2} + \frac{2m}{y^2 + 1}
$$
where the last equality holds due to the relation between $y$ and $m$.  The cases
$m < -1$ or $m > 1$ can be reduced to this setup through a rescaling of
the circulations and a relabeling of the vortices.

\begin{figure}[tb]
\centering
\includegraphics[height=8cm,keepaspectratio=true]{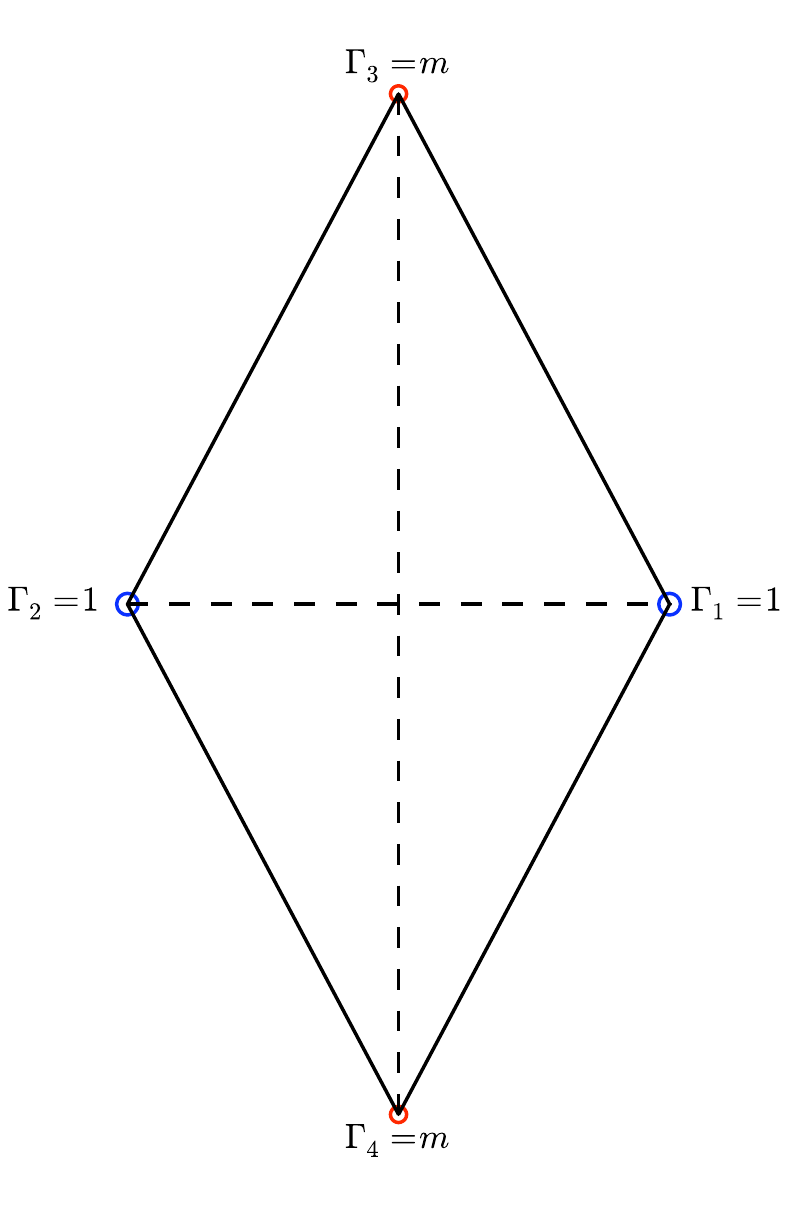}
\hspace{0.3in}
\includegraphics[height=3cm,keepaspectratio=true]{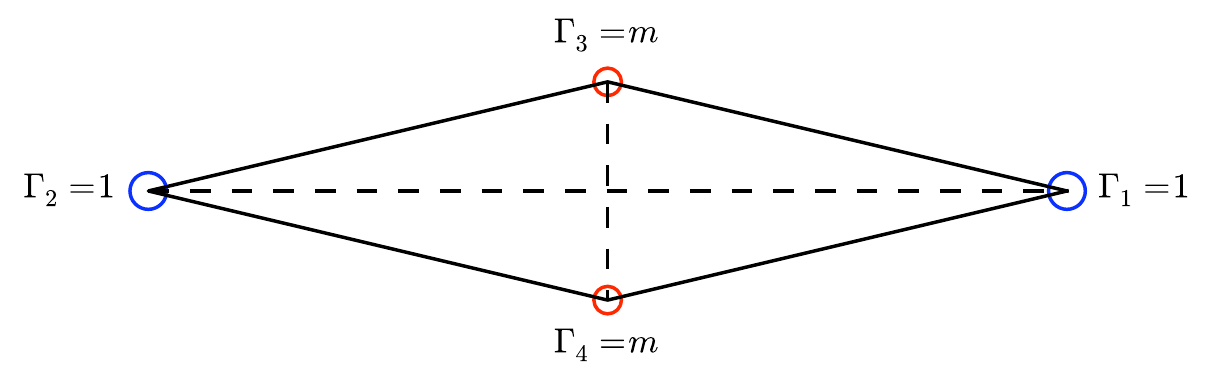}
\caption{The two distinct rhombi relative equilibria when $m=-0.2$.  The solutions rotate in opposite directions.
The configuration on the left (A) is stable; the one on the right (B) is unstable.}
\label{fig:rhombi} 
\end{figure}

There are two solutions depending on the sign choice for $y^2$ (see Figure~\ref{fig:rhombi} for an example).
Taking $+$ in~(\ref{eq:rhom-my}) yields a solution
for $m \in (-1,1]$ that always has $\omega > 0$.  We will call this solution rhombus A.
Taking $-$ in~(\ref{eq:rhom-my}) yields a solution
for $m \in (-1,0)$ that has $\omega > 0$ for $m \in (-2 + \sqrt{3}, 0)$, but
$\omega < 0$ for $m \in (-1, -2 + \sqrt{3})$.  We will call this solution rhombus B.
At $m = -2 + \sqrt{3}$, rhombus B actually
becomes an equilibrium point ($\omega = 0$ or $\nabla H(z_0) = 0$), while rhombus A is a degenerate relative equilibrium (since $L = 0$ and
$\omega \neq 0$, Theorem~\ref{thm:degen} applies).

%----------------------------------------------------
%------------------------------------------------

\begin{figure}[bth]
\centering
\includegraphics[height=7cm,keepaspectratio=true]{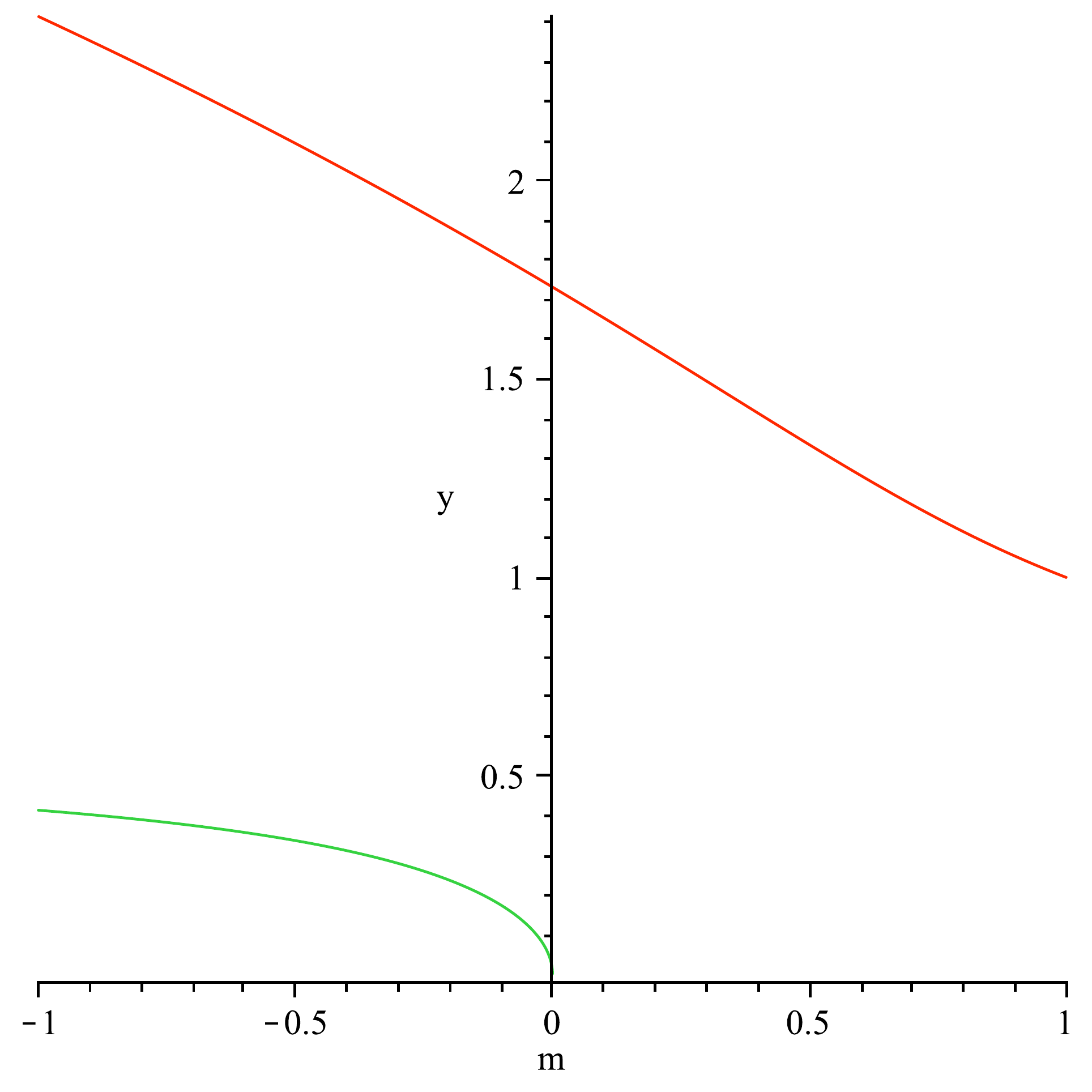}
\caption{The value of the position variable $y$ as a function of $m$. The upper curve corresponds to rhombus A
while the lower curve represents rhombus B.}
\label{fig:rhom-yvm} 
\end{figure}

It is straight forward to check that for either sign, $y^2$ is a monotonically decreasing function of~$m$.
For rhombus A, $y^2$ has a supremum of $3 + 2 \sqrt{2}$ at $m = -1$ and a minimum of $1$ at $m=1$ (the
square configuration).  
For rhombus B, $y^2$ has a supremum of $3 - 2 \sqrt{2}$ at $m = -1$ and an infimum of $0$ at $m = 0$
(vortices 3 and~4 collide).
At the bifurcation value $m = -2 + \sqrt{3}$, we have $y^2 = 2 + \sqrt{3}$ for rhombus A and
$y^2 = 7 - 4 \sqrt{3}$ for rhombus B.  The graph of $y$ as a function of $m$ for each family is shown in
Figure~\ref{fig:rhom-yvm}.

%----------------------------------------------------
%------------------------------------------------

As explained in~\cite{HRS}, the rhombus B family undergoes a pitchfork bifurcation 
at $m = m^\ast \approx -0.5951$, where $m^\ast$ is the only real root of the cubic $9m^3 + 3m^2 + 7m + 5$.  
As $m$ increases through $m^\ast$, rhombus B bifurcates into
two convex kite configurations.  We show here that one pair of eigenvalues changes from
pure imaginary to real as $m$ passes through $m^\ast$, as expected at such a bifurcation.

The symmetry of the configuration allows for an exact computation of the eigenvectors and eigenvalues
of $M^{-1} D^2H(z_0)$.   The calculations were done by hand and confirmed using Maple~\cite{maple}.
The $2 \times 2$ off-diagonal blocks of $D^2H(z_0)$ are given by
$$
A_{12} \; = \; \frac{1}{4}
\begin{bmatrix}
-1 & 0 \\
0 & 1
\end{bmatrix}, \qquad
A_{34} \; = \; \frac{m^2}{4y^2}
\begin{bmatrix}
1 & 0 \\
0 & -1
\end{bmatrix},
$$
$$
A_{13} \; = \;  \frac{m}{(y^2 + 1)^2}
\begin{bmatrix}
y^2 - 1 & 2y \\
2y & 1- y^2
\end{bmatrix}, \qquad 
A_{14} \; = \;  \frac{m}{(y^2 + 1)^2}
\begin{bmatrix}
y^2 - 1 & -2y \\
-2y & 1- y^2
\end{bmatrix},
$$
$A_{23} = A_{14}$ and $A_{24} = A_{13}$.
Due to this nice structure, it is possible to determine the nontrivial eigenvectors
of $M^{-1} D^2H(z_0)$.  They are real for all $m$ and given by
\begin{eqnarray*}
v_1 & = &  [my, 0, -my, 0, \, 0, -1, 0 , 1]^T ,\\
v_2 & = & [m, 0, m, 0, \, -1, 0 , -1 , 0]^T ,
\end{eqnarray*}
$Kv_1$ and $K v_2$.  One can check that these four vectors are $M$-orthogonal to each other and to the trivial
eigenvectors.  For $m > 0$, the vector $v_1$ corresponds to a perturbation that keeps the 
symmetry of the configuration but moves vortices $1$ and~$2$ away from the origin while pushing vortices
$3$ and~$4$ toward the origin.   This is precisely the kind of perturbation that causes instability
for the $1+n$-gon relative equilibrium (particularly when $n$ is even) in the Newtonian problem~\cite{g:stab-ring}.

The corresponding eigenvalues for $v_1$ and $v_2$, respectively,  are 
\begin{eqnarray*}
\mu_1 & = &  \frac{1}{2} - \frac{2}{(y^2 + 1)^2} \left( m(y^2 - 1) + 2 \right) \; = \;    \frac{ 7y^4 - 18y^2 + 7}{2(y^2 + 1)(3y^2 - 1)} , \\[0.1in]
\mu_2 & = &  \frac{2(m+1)(1 - y^2)}{ (y^2 + 1)^2 }  \; = \;  \frac{2(y^2-1)(y^2+2y-1)(y^2-2y-1)}{(y^2+1)^2 (3y^2-1)} .
\end{eqnarray*}
We note that the second expressions for $\mu_1$ and $\mu_2$ are well-defined 
since $y^2 \geq 1 > 1/3$ for rhombus family A and $y^2 < 3 - 2\sqrt{2} < 1/3$ for rhombus family B.
Since we have a full set of $M$-orthogonal eigenvectors for $M^{-1} D^2H(z_0)$, we can apply Lemma~\ref{lemma:charpoly2}
to obtain the four nontrivial eigenvalues.  Using equation~(\ref{eq:rhom-myA}), the key quantities for determining linear stability are
\begin{eqnarray*}
\omega^2 - \mu_1^2 & = &  \frac{ - 4(y^4 - 4y^2 + 1)(3y^4 - 2y^2 + 3)}{(y^2 +1)^2 (3y^2 - 1)^2  }   ,  \\[0.1in]
\omega^2 - \mu_2^2 & = &  \frac{y^6 \, g(y^2) \, g(1/y^2) }{4(y^2 +1)^4 (3y^2 - 1)^2  }  , \quad  g(y) = 3y^3 - 15y^2 + 41y - 5.
\end{eqnarray*}
Both of these quantities need to be positive to have linear stability.

\begin{theorem}
\begin{enumerate}
\item  Rhombus A is linearly stable for $-2 + \sqrt{3}  <  m \leq 1$.  At $m = -2 + \sqrt{3}$, the relative equilibrium
is degenerate.  For $-1 < m < -2 + \sqrt{3}$, rhombus A is unstable and the nontrivial eigenvalues consist of a real pair and a pure imaginary pair.

\item  Rhombus B is always unstable.   One pair of eigenvalues is always real.  The other pair of eigenvalues
is pure imaginary for $-1 < m < m^\ast$ and real for $m^\ast < m < 0$, where
$m^\ast \approx -0.5951$ is the only real root of the cubic $9m^3 + 3m^2 + 7m + 5$.  
At $m = m^\ast$, rhombus B is degenerate.
\end{enumerate}
\end{theorem}

\pf
Consider $\omega^2 - \mu_1^2$ as a function of $y$.  Since $3y^4 - 2y^2 + 3 = 3(y^2 - 1/3)^2 + 8/3 > 0$,
the sign of $\omega^2 - \mu_1^2$ is determined by the sign of $y^4 - 4y^2 + 1$.  Note that this quartic
has roots when $y^2 = 2 \pm \sqrt{3}$.  It follows that the pair of eigenvalues corresponding to $\{v_1, Kv_1 \}$
are pure imaginary if $2 - \sqrt{3} < y^2 < 2 + \sqrt{3}$, real if $0 < y^2 < 2 - \sqrt{3}$ or $y^2 > 2 + \sqrt{3}$,
and $0$ when $y^2 = 2 \pm \sqrt{3}$.  Since $y^2$ is monotonically decreasing for both families, we can easily translate
these statements in terms of~$m$. For rhombus A, we have that $1 \leq y^2 < 2 + \sqrt{3}$ for
$1 \geq  m > -2 + \sqrt{3}$, so the eigenvalues are pure imaginary in this interval of $m$-values.
When $-2 + \sqrt{3} > m > -1$, $y^2 > 2 + \sqrt{3}$, and the eigenvalues form a real pair.  At $m = -2 + \sqrt{3}$,
the vector $Kv_1$ is in the kernel of the stability matrix $B$ so rhombus A becomes degenerate, which is in agreement with
Theorem~\ref{thm:degen} since $L = 0$ here.  For rhombus B, we have $y^2 < 3 - 2 \sqrt{2} < 2 - \sqrt{3}$
for $-1 < m < 0$, so the eigenvalues always form a real pair.

Next we consider the eigenvalues associated to $\{v_2, Kv_2 \}$.  In this case, the sign of $\omega^2 - \mu_2^2$ is 
determined by the sign of the product $g(y^2)  g(1/y^2)$.  The discriminant of the cubic $g$ is negative, so
it has only one real root which we denote by $\kappa$.  
Computing a Gr\"{o}bner basis for $g(y^2)$ and equation~(\ref{eq:rhom-myA}) yields
$\kappa = -m^\ast (3 m^\ast + 2) \approx 0.1278$.  
Moreover, substituting $y^2 = \kappa$ into equation~(\ref{eq:rhom-myA}) gives $m = m^\ast$.    
The pair of eigenvalues corresponding to $\{v_2, Kv_2 \}$ are pure imaginary for $\kappa < y^2 < 1/\kappa$, 
real for $0 < y^2 < \kappa$ or $y^2 > \kappa$, and $0$ when $y^2 = \kappa$ or $1/\kappa$.
For rhombus A, we have $\kappa < 1 \leq y^2 <  3 + 2 \sqrt{2} < 1/\kappa$, where the last inequality is rigorously
shown by checking $g(1/(3 + 2\sqrt{2})) > 0$.  Thus, this pair of eigenvalues is always pure imaginary for the rhombus A
family.   For rhombus B, we have $\kappa < y^2 < 1/3 < 1/\kappa$ as long as $-1 < m < m^\ast$.  For these $m$-values,
the pair of eigenvalues is pure imaginary.  At $m = m^\ast$, $y^2 = \kappa$ and rhombus B becomes a degenerate relative equilibrium
with $v_2$ in the kernel of the stability matrix $B$.
For $m^\ast < m < 0$, $y^2 < \kappa$ and the eigenvalues bifurcate into a real pair.  Combining these conclusions with
those for the eigenvectors $\{v_1, Kv_1 \}$ finishes the proof.
\enpf

\begin{remark}
\begin{enumerate}
\item
For completeness, we note that the equilibrium point for rhombus B at $m = -2 + \sqrt{3}$ is linearly unstable.
This follows since the nontrivial eigenvalues of $M^{-1}D^2H(z_0)$ are the two real pairs
$\pm \mu_1$ and $\pm \mu_2$.
\item  
If $m > 0$, then Theorem~\ref{thm:stable} implies that the rhombus A family is actually nonlinearly stable.
This provides a possible explanation for the long-term persistence of the rhombus-like configuration
in the numerical experiments of Kossin and Schubert~\cite{KossSchub}.
\item
For $-2 + \sqrt{3} < m < 0$, rhombus A is linearly stable, but numerical calculations using Matlab show that
it is a saddle of $H$ restricted to $I = I_0 > 0$.  This shows that for opposite signed circulations, 
it is possible for a relative equilibrium to be linearly stable, but not a minimum of $H$ restricted to a level surface
of~$I$.  The proof of Theorem~\ref{thm:main} breaks down because $v_1^T M v_1$ and $v_2^T M v_2$ are
each negative, so the value of $q_{z_0}(v)$ can also be negative, leading to a saddle.
\end{enumerate}
\end{remark}

%%%%%%%%%%%%%%%%%%%%%%%%%%%%%%%%%%%%%%%%%
\subsection{The Isosceles Trapezoid Family}
%%%%%%%%%%%%%%%%%%%%%%%%%%%%%%%%%%%%%%%%%

Another one-parameter family of symmetric relative equilibria in the four-vortex problem is
formed by isosceles trapezoid configurations.
This family was described and analyzed in~\cite{HRS}.  
Set  $\Gamma_1 = \Gamma_2 = 1$ and $\Gamma_3 = \Gamma_4 = m$, treating $m \in \mathbb{R}^+$ as a parameter.
Vortices 1 and~2 lie on one base of the trapezoid while vortices 3 and~4 lie on the other, and the two legs
of the trapezoid are congruent.   The pair of vortices with the larger circulations lie on the longer base, with $m=1$ corresponding to a square.
This family exists only for $m > 0$.  As $m \rightarrow 0$, the configuration approaches an equilateral triangle with
vortices 3 and 4 colliding in the limit.

\begin{figure}[htb]
\centering
\includegraphics[height=7cm,keepaspectratio=true]{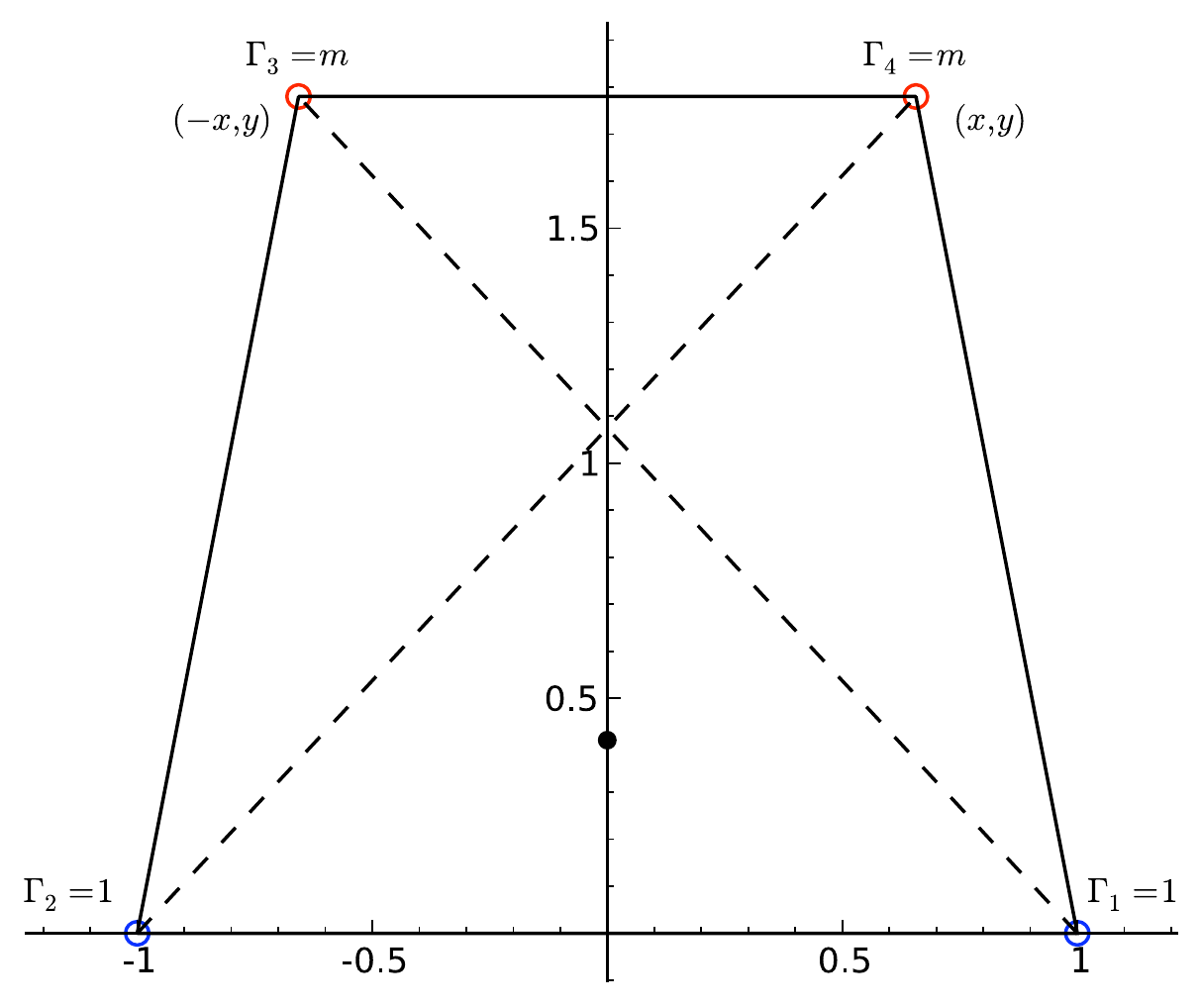}
\caption{The setup for the isosceles trapezoid family.  Figure shown is for $m = 0.3$, with the center
of vorticity indicated by the small solid disk.}
\label{fig:trap-fig} 
\end{figure}

The coordinates of the trapezoid are given by $\hat{z_1} = (1,0), \hat{z_2} = (-1,0),$
$\hat{z_3} = (-x,y)$ and  $\hat{z_4} = (x,y)$, which places the center of vorticity at $\hat{c} = (0, my/(m+1))$ (see Figure~\ref{fig:trap-fig}).
To be a relative equilibrium, we must have 
$$
(x, y) \; = \; ( \sqrt{\alpha}, \sqrt{2m + 3 - \alpha}) \quad \mbox{where } \quad \alpha =  \frac{m(m+2)}{2m+1}.
$$
These expressions were derived from the corresponding results in~\cite{HRS}.  One can check that
the positions $z_i = \hat{z_i} - c$ satisfy the defining equations~(\ref{eq:rel-equ}) with angular velocity $\omega = (2m+1)/(2m+2)$.
This value of $\omega$ concurs with the one obtained using the formula $\omega = L/(2I)$.  Note that $x$ and $y$ are both real
only when $m > 0$.

%----------------------------------------------------
%------------------------------------------------

Since the key matrix $M^{-1} D^2H(z_0)$ is invariant under translation, we can evaluate it by substituting in the
values $\hat{z_i}$ given above.  Due to symmetry, it is possible to give explicit formulas for the eigenvectors of $M^{-1} D^2H(z_0)$.
Note that $r_{13}^2 = r_{24}^2 = (x+1)^2 + y^2$ and $r_{14}^2 = r_{23}^2 = (1-x)^2 + y^2$.
The $2 \times 2$ off-diagonal blocks of $D^2H(z_0)$ are given by
$$
A_{12} =  \frac{1}{4}
\begin{bmatrix}
-1 & 0 \\
0 & 1
\end{bmatrix}, \qquad
A_{34} \; = \; \frac{m^2}{4x^2}
\begin{bmatrix}
-1 & 0 \\
0 & 1
\end{bmatrix},
$$
$$
A_{13} =   \frac{m}{r_{13}^4}
\begin{bmatrix}
y^2 - (x+1)^2 & 2y(x+1) \\
2y(x+1) & (x+1)^2 - y^2
\end{bmatrix}, \qquad 
A_{24} \; = \;  \frac{m}{r_{13}^4}
\begin{bmatrix}
y^2 - (x+1)^2 & -2y(x+1) \\
-2y(x+1) & (x+1)^2 - y^2
\end{bmatrix},
$$

$$
A_{14} =  \frac{m}{r_{14}^4}
\begin{bmatrix}
y^2 - (1-x)^2 & 2y(1-x) \\
2y(1-x) & (1-x)^2 - y^2
\end{bmatrix}, \qquad
A_{23} \; = \;  \frac{m}{r_{14}^4}
\begin{bmatrix}
y^2 - (1-x)^2 & -2y(1-x) \\
-2y(1-x) & (1-x)^2 - y^2
\end{bmatrix}.
$$

Surprisingly, the null space of $D^2H(z_0)$ is four-dimensional for all values of $m$.  In addition to the 
two vectors coming from the conservation of the center of vorticity, the vectors $w_1, Kw_1$ are also
in the kernel of $D^2H(z_0)$, where $w_1 =   [mx, 0, -mx, 0, \, 1, 0, -1 , 0]^T $.  The vector $w_1$ was found
by searching for vectors in the kernel of the form $[a, 0, -a, 0, \, b, 0, -b , 0]^T$, and was verified symbolically using Maple.
Taking $\mu_1 = 0$ in Lemma~\ref{lemma:charpoly2} shows that two of the nontrivial eigenvalues 
for the isosceles trapezoid family are $\pm \omega  i$.  Note that the invariant space formed by the span
of $w_1$ and $Kw_1$ is $M$-orthogonal to the invariant space arising from the conservation of the center of vorticity.
Thus, although the eigenvalues $\pm \omega  i$ are repeated, the Jordan form of $B$ does not have an off-diagonal block.

The remaining nontrivial eigenvectors for $M^{-1} D^2H(z_0)$ are harder to locate.  After some initial numerical investigations using Maple,
we found the pair of eigenvectors $w_2, Kw_2$ with corresponding eigenvalues $\pm \mu_2$, 
where $w_2 = [-a, -m, a, -m, \, ax, 1, -ax , 1]^T $, 
$$
a \; = \;  \frac{m \sqrt{ 3(2m+1)} }{m^2 + m + 1} 
\quad \mbox{ and }  \quad
\mu_2 \; = \;  \frac{ m^2+m+1}{2(m+1)(m+2)}.
$$
This, in turn, gives
$$
\omega^2 - \mu^2  \; = \;  \frac{3(m^2 + 4m + 1)}{4(m+2)^2}  \; = \;   \frac{3L}{4(m+2)^2}  \; > \; 0.
$$
It is interesting, though not altogether surprising, to see the reappearance of the total angular vortex momentum~$L$ here.
It follows that the remaining nontrivial eigenvalues lie on the imaginary axis and the isosceles trapezoid family is stable for
all~$m$.  We have proven the following result.
\begin{theorem}
The isosceles trapezoid family is stable for all values of the parameter $m > 0$.  The nontrivial eigenvalues are
$\pm \omega  i$ and $\displaystyle{ \frac{\pm \sqrt{3L}}{2(m+2)} \, i }$.
\end{theorem}

%%%%%%%%%%%%%%%%%%%%%%
\section{Conclusion}
%%%%%%%%%%%%%%%%%%%%%%

We have adapted the approach of Moeckel from the $n$-body problem to study the linear stability 
of relative equilibria in the planar $n$-vortex problem.  
Due to some special properties of the logarithmic Hamiltonian in the vortex case,  
the stability is easier to determine and a complete factorization of the characteristic
polynomial into quadratic factors exists when all the circulations have the same sign.  
Using a topological approach motivated by a conjecture of Moeckel's, we show that,
in the case of same-signed circulations, a relative equilibrium is linearly stable if and only
if it is a nondegenerate minimum of the Hamiltonian restricted to a level surface of the
angular impulse.  In this case, it follows quickly that linear stability implies nonlinear stability as
well.  Two symmetric examples of stable configurations in the four-vortex problem, a rhombus and an
isosceles trapezoid, are analyzed in detail, with interesting bifurcations discussed.
It is hoped that this work will be of some use to researchers in fluid dynamics, particularly
to those who have discovered vortex crystals in their numerical simulations.

\vs

\noindent  {\bf Acknowledgments:} 
The author would like to thank the National Science Foundation (grant DMS-1211675) and 
the Banff International Research Station for their support, as well as Anna Barry, Florin Diacu
and Dick Hall for stimulating discussions related to this work.

%%%%%%%%%%%%%%%%%%%%%%%%%%%%%%%%%%%%%%%%%%%%%%%%%%%%%%%%%%%%%%%%%%%%%%%%%%%%%%%%%%
\bibliographystyle{amsplain}

\end{document}